\tolerance=10000
\raggedbottom

\baselineskip=15pt
\parskip=1\jot

\def\sk{\vskip 3\jot}

\def\heading#1{\vskip3\jot{\noindent\bf #1}}
\def\label#1{{\noindent\it #1}}


\def\ref#1;#2;#3;#4;#5.{\item{[#1]} #2,#3,{\it #4},#5.}
\def\refinbook#1;#2;#3;#4;#5;#6.{\item{[#1]} #2, #3, #4, {\it #5},#6.} 
\def\refbook#1;#2;#3;#4.{\item{[#1]} #2,{\it #3},#4.}


\def\({\bigl(}
\def\){\bigr)}

\def\ga{\gamma}
\def\ze{\zeta}
\def\la{\lambda}
\def\si{\sigma}
\def\ps{\psi}

\def\Ga{\Gamma}

\def\calP{{\cal P}}
\def\calQ{{\cal Q}}
\def\calR{{\cal R}}

\def\calW{{\cal W}}

\input xy
\xyoption{all}

\def\edge#1{\ar@{-}[#1]}

\def\text#1{{\rm #1}}
\def\state#1{*+++[o][F]{#1}}
\def\fstate#1{*++++[o][F=]{#1}}

\def\Ex{{\rm Ex}}
\def\Var{{\rm Var}}

\def\Stwo#1#2{\left\{#1\atop#2\right\}}

\def\fcfs{{\rm FCFS}}
\def\lcfs{{\rm LCFS}}
\def\ros{{\rm ROS}}
\def\fos{{\rm FOS}}

{
\pageno=0
\nopagenumbers
\rightline{\tt egp.ex.tex}
\vskip1in

\centerline{\bf Analysis of an $M/M/1$ Queue Using Fixed Order of Search for Arrivals and Service}
\vskip0.5in

\centerline{Patrick Eschenfeldt}
\centerline{\tt peschenfeldt@hmc.edu}
\sk

\centerline{Ben Gross}
\centerline{\tt bgross@hmc.edu}
\sk

\centerline{Nicholas Pippenger}
\centerline{\tt njp@math.hmc.edu}
\sk

\centerline{Department of Mathematics}
\centerline{Harvey Mudd College}
\centerline{1250 Dartmouth Avenue}
\centerline{Claremont, CA 91711}
\vskip0.5in

\noindent{\bf Abstract:}
We analyze an $M/M/1$ queue with a service discipline in which customers, upon arriving when the server is busy, search a sequence of stations for a vacant station at which to wait, and in which the server, upon becoming free when one or more customers are waiting, searches the stations in the same order for a station occupied by a customer to serve.
We show how to find complete asymptotic expansions for all the moments of the waiting time
in the heavy traffic limit.
We show in particular that the variance of the waiting time for this discipline is more similar to that of last-come-first-served (which has a pole of order three as the arrival rate approaches the service rate) than that of first-come-first-served (which has pole of order two).
\vskip0.5in
\leftline{{\bf Keywords:} Queueing theory,  Lambert series, asymptotic expansions.}
\sk
\leftline{{\bf Subject Classification:} 60K26, 90B22}

\vfill\eject
}

\heading{1. Introduction}

We consider the $M/M/1$ queue (with independent exponentially distributed interarival times,
independent exponentially distributed service times, and a single server) with various service disciplines.
We shall be interested mainly in the ``heavy traffic'' limit, $\la\to 1$, where $\la$ is the arrival rate (measured in units of the service rate), and all asymptotic statements in this paper refer to this limit.

It is well known (see Little [L]) that the average waiting time $\Ex[W]$ ($W$ is the length of interval from arrival to commencement of service) does not depend on the service discipline (the rule used to determine which waiting customer is served next when the server becomes free).
We have
$$\Ex[W] = {\la\over 1-\la}\sim {1\over 1-\la}. \eqno(1.1)$$

The variance of $W$, however, (and more generally its higher moments) does depend on the service discipline.
Kingman [K] has shown that, among all service disciplines, ``first-come-first-served'' (FCFS) minimizes the variance of $W$.
We have
$$\Var[W_\fcfs] = {\la(2-\la) \over (1-\la)^2} \sim {1\over (1-\la)^2} \eqno(1.2)$$
(see for example Riordan [R2, pp.102--103]).
Tambouratzis [T] has shown that ``last-come-first-served'' (LCFS) maximizes this variance.
We have
$$\Var[W_\lcfs] = {\la(2-\la+\la^2)\over (1-\la)^3}\sim {2\over (1-\la)^3} \eqno(1.3)$$
(see for example Riordan [R2, pp.~106--109]; LCFS was first analyzed by Vaulot [V2]).
We note that the difference between FCFS and LCFS is qualitative, in that $\Var[W_\lcfs]$ has a pole of order three at $\la=1$, whereas $\Var[W_\fcfs]$ has a pole of order only two there.

Another service discipline that has been studied is ``random-order-of-service'' (ROS),
first successfully analyzed by Vaulot [V1].
One of the motivations for studying ROS was stated by Riordan [R1]:
{
\medskip\narrower
``In many switching systems it is not feasible to fully realize this ethical ideal of first come, first served, and it has long been of interest to determine delays on another basis.
The contrasting assumption is of calls picked at random, which is again an idealization but in large offices appears to be called for, as a bound for the service actually given.''
\medskip
}
\noindent
This statement suggests that (1) in practical systems of that era (around 1953) it was not possible to keep track of the order of arrival, (2) ROS was analyzed as a substitute for the service discipline actually implemented, and (3) it was hoped that the performance of ROS would approximate that of the discipline actually implemented.
For ROS, we have
$$\Var[W_\ros] = {\la(4-2\la+\la^2) \over (2-\la) \, (1-\la)^2} \sim {3\over (1-\la)^2} \eqno(1.4)$$
(see for example Riordan [R2, pp.~103--106]).
Comparing (1.4) with (1.2) and (1.3), we see that ROS is 	qualitatively ``more like'' FCFS than LCFS, in that $\Var[W_\ros]$ has a pole of order two rather than three, though its coefficient is larger than that of FCFS by a factor of three.

In the early 1950s, telephone switching systems were electromechanical and did not employ randomization beyond that present in the arrival and service processes.
In this paper we shall analyze a service discipline that we shall call ``fixed-order-of-search'' (FOS).
For this discipline, there is an infinite sequence $\calW_1, \calW_2, \ldots$ of ``waiting stations'',
each of which can be either ``vacant'' or ``occupied''.
A customer arriving when the server is busy searches these stations in increasing order of their indices and occupies the first vacant station it finds.
When the server becomes free and one or more customers are waiting, it searches the stations in the same order and serves the customer waiting at the first occupied station it finds, thereby vacating that station.
This discipline was introduced by Eschenfeldt, Gross and Pippenger [E].
Like FCFS and LCFS (and unlike ROS), it does not employ randomization beyond that present in the arrival and service processes.
If the server were to search the stations in the reverse order to that used by arriving customers
(serving the customer waiting at the occupied station with the largest, rather than the smallest, index), the result would be LCFS, with its concomitant maximum variance for the waiting time.
The choice of the same order of search for both customers and the server thus represents an attempt to improve upon LCFS, while still using a fixed order of search in each case.

We shall give an exact formula for $\Var[W_\fos]$: 
$$\Var[W_\fos] = 
{\la\,(6+\la+\la^2) \over (1-\la)^3} - {4\,\ps^{(1)}_\la(1) \over (1-\la)\,\log^2 \la}, \eqno(1.5)$$
where $\ps^{(1)}_q(x)$ is the $q$-trigamma function (defined below).
We shall indicate how similar, but increasingly more complicated, formulas can be derived for the higher moments of $W_\fos$.
We shall also give an asymptotic formula for $\Var[W_\fos]$:
$$\Var[W_\fos] \sim  {8 - 4\,\ze(2) \over (1-\la)^3} \eqno(1.6)$$
where $\ze(s) = \sum_{n\ge 1} 1/n^s$ is the Riemann zeta function
and $\ze(2)=\pi^2/6$
(see for example Whittaker and Watson [W2, pp.~265--280]).
We shall also indicate how complete asymptotic expansions (with error terms of the form 
$O\((1-\la)^R\)$ for any $R$) can be derived for the variance of $W_\fos$, as well as for the higher moments.
Comparing (1.6) with (1.2) and (1.3), we see that FOS is 	qualitatively ``more like'' LCFS than FCFS, in that $\Var[W_\fos]$ has a pole of order three rather than two, though its coefficient is smaller than that of LCFS by a factor of $4-2\,\ze(2) = 4-\pi^2/3 = 0.7101\ldots\,$.

Eschenfeldt, Gross and Pippenger [E] initiated the study of FOS, determining the distribution of the index $I$ of the station $\calW_I$ at which a newly arriving customer waits (where $I=0$ if the server is idle at the time of the arrival):
we have $\Pr[I\ge 0] = 1$ and
$$\Pr[I\ge i] = {(1-\la) \, \la^{i} \over 1-\la^{i}} \eqno(1.7)$$
for $i\ge 1$.
The moments of $I$ can be expressed in terms of the sums
$$T_l(\la) = \sum_{i\ge 1} {i^l \, \la^j \over 1 - \la^i}, \eqno(1.8)$$
for which we have the asymptotic formulas
$$T_0(\la)  \sim {1\over 1-\la}\log {1\over 1-\la} \eqno(1.9)$$
and
$$T_l(\la)  \sim {l! \, \ze(l+1)\over (1-\la)^{l+1}} \eqno(1.10)$$
for $l\ge 1$, and the exact formulas
$$T_0(\la) = {\ps_\la(1) + \log(1-\la) \over \log\la} \eqno(1.11)$$
and
$$T_l(\la) = {\ps_\la^{(l)}(1) \over \log^{l+1} \la} \eqno(1.12)$$
for $l\ge 1$.
Here $\ps_q(x) = \partial\log\Ga_q(x)/\partial x$ is the $q$-digamma function, the logarithmic derivative of the $q$-gamma function 
$\Ga_q(x) = (1-q)^{1-x} \,\prod_{n\ge 0} \((1-q^{n+1})/(1-q^{n+x})\)$
(see for example Gasper and Rahman [G, p.~16]),
and $\ps_q^{(l)}(x) = \partial^l \ps_q(x)/\partial x^l$ is the $l$-th $q$-polygamma function.
(The sums $T_l(\la)$, which are called Lambert series (see for example Hardy and Wright [H, p.~257]), are the generating functions
$T_l(\la) = \sum_{n\ge 1} \si_l(n)\,\la^n$ for the sums $\si_l(n) = \sum_{d\mid n} d^l$
of the $l$-th powers of the divisors of of $n$ (see for example Hardy and Wright [H, p.~239]).)
In terms of the $T_l(\la)$, we have
$$\Ex[I^m] = (1-\la) \sum_{0\le l\le m-1} {m\choose l} (-1)^{m-1-l} \, T_l(\la).$$

In Section 2 we shall determine the moment generating $M_W(s)$ function for $W_\fos$ (which  in what follows we shall denote simply $W$).
In Section 3, we shall derive the
exact formula (1.5)
and the asymptotic formula (1.6).
We shall also indicate how similar exact and asymptotic formulas can be found for the higher moments of  $W$.
Finally, we shall indicate how these asymptotic formulas can be extended to complete asymptotic expansions (with error terms of the form $O\((1-\la)^R\)$ for any $R$) for these quantities.
\sk

\heading{2. The Generating Functions}

Consider the random process whose state variable $J$ denotes the number of customers in the system.
The random variable $J$ is zero during an idle period (interval of time during which the server is idle).
It is incremented whenever a customer arrives, and decremented whenever a customer departs
(that is, at the termination of a service interval).
Arrivals occur in a Poisson process with rate $\la$.
During a busy period (interval of time during which $J\ge 1$), departures occur 
in an independent Poisson process at rate $1$.
Thus, during a busy period,  ``transitions'', by which we mean arrivals and departures together,
occur in a Poisson process at rate $1+\la$.
Furthermore, during a busy period, the probability that the next transition will be an arrival is 
$p = \la/(1+\la)$, and the probability that it will be a departure is $q = 1/(1+\la)$.
In this section we shall study the distribution of the random variable $N$, defined as the number of transitions that occur between the arrival of a customer (excluded) and the departure that initiates its service interval (included).
Specifically, we shall determine the probability generating function 
$g(t) = \sum_{n\ge 0} \Pr[N = n] \, t^n$
for $N$.
We have $N=0$ if the arrival initiates a busy period, and $N\ge 1$ if it occurs during a busy period.

The index $I$ of the station $\calW_I$ at which a newly arriving customer waits has the distribution given by (1.7).
As a first step to determining $g(t)$, we shall 
determine the conditional generating function 
$g_i(t) = \sum_{n\ge 0} \Pr[N = n\mid I = i] \, t^n$.
We begin with two special cases.
If $i=0$, the arrival initiates a busy period, so $N=0$ and
$g_0(t) = 1$.
If $i=1$, the customer waits at $\calW_1$ and will be served as soon as the next departure occurs.
The number of transitions preceding and including this next departure has a geometric  distribution, and $g_1(t)  = \sum_{n\ge 1} p^{n-1} \, q \, t^n = qt/(1-pt)$.

For the general case, we consider the random process whose state variable $K$ denotes the number of vacant stations among $\calW_1, \ldots, \calW_i$.
Since the customer in question waits at station $\calW_i$, we have $K=0$ immediately after the arrival of that customer.
Furthermore, the first time thereafter at which $K=i$ coincides with the beginning of the service interval for that customer, and thus occurs after exactly $N$ transitions have occurred.
When $1\le K\le i-1$, the random variable $K$ is incremented by a departure (because the next customer served will be waiting at one of the stations under consideration) and decremented by an arrival (because the arriving customer will wait at one of these stations).
If however $K=0$, a departure will increment $K$, but
an arrival will leave $K$ unchanged (because it will have to wait at a station beyond $\calW_i$).
Thus the process determining $N$ given $I=i$ is a random walk with one ``reflecting barrier''
(at $K=0$) and one ``absorbing barrier'' (at $K=i$), as shown in Figure 1.
\medskip
$$\xymatrix{
\ar[r]^-{\text{start}} & \state{0} \ar@(ul,ur)[]^p \ar@/^/[r]^q & \state{1} \ar@/^/[r]^-q \ar@/^/[l]^p&\cdots \ar@/^/[r]^-q \ar@/^/[l]^-p & \state{\scriptstyle i-1} \ar[r]^q \ar@/^/[l]^-p & \fstate{\scriptstyle \text{stop}}
}$$
\medskip
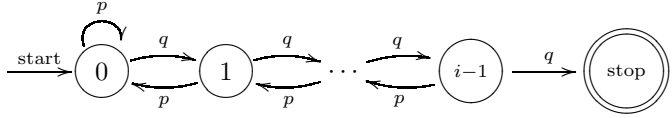
\centerline{Figure 1. The process determining $N$, given that $I=i$.} 
\centerline{Non-terminal states are labeled with values of the random variable $K$.}
\medskip
The number $N$ of steps to absorption in this process has been studied by Weesakul [W1],
who shows that the generating function is given by
$$g_i(t) = {q^i \, t^i \, \(Q(t) - P(t)\) \over Q(t)^{i+1} - P(t)^{i+1} - pt\(Q(t)^i - P(t)^i\)}, \eqno(2.1)$$
where
$$P(t) = {1 - \sqrt{1 - 4pqt^2} \over 2}$$
and
$$Q(t) = {1 + \sqrt{1 - 4pqt^2} \over 2}.$$
We note that for $t=1$ we have $P(1) = p$, $Q(1) = q$ and $g_i(1) = 1$.

We can now express the unconditional generating function $g(t)$ by using summation by parts:
$$\eqalign{
g(t) 
&= \sum_{i\ge 0} g_i(t) \, \Pr[I = i] \cr
&= \sum_{i\ge 0} g_i(t) \, (\Pr[I \ge i] - \Pr[I \ge i+1]) \cr
&= g_0(t) \, \Pr[I\ge 0] + \sum_{i\ge 1} \(g_i(t) - g_{i-1}(t)\) \, \Pr[I \ge i] \cr
&= 1 + \sum_{i\ge 1} \(g_i(t) - g_{i-1}(t)\) \, \Pr[I \ge i], \cr
}$$
because $g_0(t) = \Pr[I\ge 0] = 1$.
Thus, using (1.7), we have
$$g(t) = 1 + (1-\la) \sum_{i\ge 1} {(\nabla g)_i(t) \, \la^i \over 1-\la^i}, \eqno(2.2)$$
where $(\nabla g)_i(t) = g_i(t) - g_{i-1}(t)$ denotes the backward difference of $g_i(t)$.
We note that for $t=1$ we have $(\nabla g)_i(1) = 0$, so $g(1) = 1$.

We are now ready to drive the moment generating function $M_W(s) = \Ex[e^{sW}]$ for $W$.
Each intertransition time $X$ is exponentially distributed with mean $1/(1+\la)$,
so the moment generating function for $X$ is $M_X(s) = (1+\la)/(1 + \la - s)$.
The moment generating function for the sum $\sum_{1\le k\le n} X_k$ of $n$ independent intertransition times $X_1,\ldots, X_n$ is $M_X(s)^n = \((1+\la)/(1 + \la - s)\)^n$.
Thus the waiting time $W$, which is the sum of the random number $N$ of independent intertransition times has the moment generating function
$$\eqalignno{
M_W(s) 
&= \sum_{n\ge 0} \Pr[N=n] \, M_X(s)^n \cr
&= g\(M_X(s)\) &(2.3)\cr
&= g\left({1+\la\over 1+\la- s}\right). \cr
}$$
\sk

\heading{3. The Moments}

In this section we shall derive the mean and variance for $N$ and $W$, and indicate how to derive the higher moments as well.

For the mean of $N$, we use the formula $\Ex[N] = g'(1)$.
We have 
$$g_i'(1) = {(1+\la) \, \(i(1-\la) - \la(1-\la^i)\) \over (1-\la)^2},$$
so
$$(\nabla g')_i(1) = {(1+\la) \, (1-\la^i) \over 1-\la}.$$
Thus we have
$$\eqalignno{
\Ex[N]
&= g'(1) \cr
&= (1-\la) \sum_{i\ge 1} {(\nabla g')_i(1) \, \la^i \over 1-\la^i} \cr
&= (1-\la) \sum_{i\ge 1} {(1+\la) \, (1-\la^i) \over 1-\la} \, {\la^i \over 1-\la^i} \cr
&= {(1+\la) \, \la \over 1-\la}. &(3.1) \cr
}$$
For the mean of $W$, we use the formula for the expectation of the sum of a random number $N$
of independent, identically distributed random variables $X, X_1, X_2, \ldots\,$:
$\Ex\left[\sum_{1\le k\le N} X_k\right] = \Ex[N] \, \Ex[X]$.
Since the intertransition time $X$ satisfies $\Ex[X] = 1/(1+\la)$,
we have 
$$\eqalign{
\Ex[W] 
&= \Ex[N] \, \Ex[X] \cr
&= {(1+\la) \, \la \over 1-\la} \, {1\over 1+\la} \cr
& = {\la \over 1-\la}, \cr
}$$
in accordance with (1.1).
 
 For the variance of $N$, we begin by using the formula for the second factorial moment:
 $\Ex[N(N-1)] = g''(1)$.
We have
 $$\eqalign{
 g''_i(1) 
 &= 
 { (1-\la)^2 (1+\la)^2\,i^2 -   (1-\la)(1+\la)(1-10\la-3\la^2)\,i \over (1-\la)^4} \cr
&\qquad -{ \(6 \la (1-\la)(1 + \la)^2\,i + 2 \la (1 + \la)  (1+4\la-\la^2) - 2\la^2(1+\la)^2 \,\la^i\) (1-\la^i) 
\over (1-\la)^4}, \cr
}$$
so
\def\qpolylin{2(1-\la)(1+\la)^3}
\def\qpolyconstlin{6i(1-\la)^2(1+\la)^2}
\def\mqpolyconstconst{2(1-\la)^3(1+\la)}
\def\qpoly{\(\qpolylin \, \la^i + \qpolyconstlin  - \mqpolyconstconst\)}
\def\mrpoly{4i(1-\la)^2(1+\la)^2}
$$(\nabla g'')_i(1) = {\qpoly \, (1-\la^i) - \mrpoly \over (1-\la)^4}.$$
Thus we have
$$\eqalign{
\Ex[N(N-1)]
&= {2(1+\la)^3\over (1-\la)^2} \sum_{i\ge 1} \la^{2i} 
+ {6(1+\la)^2\over 1-\la} \sum_{i\ge 1} i\,\la^i
- 2(1+\la) \sum_{i\ge 1} \la^i 
- {4(1+\la)^2\over 1-\la} \sum_{i\ge 1} {1\, \la^i \over 1-\la^i} \cr
&= {2\la^2\,(1+\la)^2 \over (1-\la)^3}
+ {6\la(1+\la)^2 \over (1-\la)^3}
- {2\la\,(1+\la) \over 1-\la}
- {4(1+\la)^2 \over 1-\la} \, T_1(\la), \cr.
}$$
where we have used the definition (1.8) to evaluate the last sum.
It follows that
$$\eqalign{
\Var[N]
&= \Ex[N(N-1)] + \Ex[N]  - \Ex[N]^2 \cr
&= {2\la^2\,(1+\la)^2 \over (1-\la)^3}
+ {6\la(1+\la)^2 \over (1-\la)^3}
- {\la\,(1+\la) \over 1-\la}
-{\la^2(1+\la)^2 \over (1-\la)^2}
- {4(1+\la)^2 \over 1-\la} \, T_1(\la), \cr.
&= {\la(1+\la)^2 \,(6+\la+\la^2) \over (1-\la)^3}
- {\la\,(1+\la) \over 1-\la}
- {4(1+\la)^2 \over 1-\la} \, T_1(\la), \cr.
}$$
where we have used (3.1), then combined the first, second and fourth terms.
For the variance of $W$, we use the formula for the variance of a random number $N$ of independent, identically distributed random variables $X, X_1, X_2, \ldots\,$:
$\Var\left[\sum_{1\le k\le N} X_k\right] = \Var[N] \, \Ex[X]^2 + \Ex[N] \, \Var[X]$.
Since the intertransition time $X$ satisfies
$\Ex[X] = 1/(1+\la)$ and $\Var[X] = 1/(1+\la)^2$, we have
$$\eqalign{
\Var[W]
&= \Var[N] \, \Ex[X]^2 + \Ex[N] \, \Var[X] \cr
&= \left({\la(1+\la)^2 \,(6+\la+\la^2) \over (1-\la)^3}
- {\la\,(1+\la) \over 1-\la}
- {4(1+\la)^2 \over 1-\la} \, T_1(\la)\right) \, {1\over (1+\la)^2}
+ {\la\,(1+\la) \over 1-\la} \, {1\over (1+\la)^2} \cr
&= {\la(6+\la+\la^2) \over (1-\la)^3} - {4\over 1-\la} \, T_1(\la).
}$$
Evaluating $T_1(\la)$ using (1.12) yields
$$\Var[W] = {\la(6+\la+\la^2) \over (1-\la)^3} - {4\,\ps^{(1)}_\la(1) \over (1-\la)\,\log^2 \la},$$
confirming (1.5), whereas using (1.10) yields
$$\Var[W]\sim {\la(6+\la+\la^2) \over (1-\la)^3} - {4\,\ze(2) \over (1-\la)^3},$$
confirming (1.6).

It is straightforward to generalize the derivations of the mean and variance of $W$ given above to the higher moments.
We begin by indicating how to derive the higher factorial moments of $N$.
We have
$$\eqalign{
\Ex[N(N-1)\cdots(N-m+1)]
&= g^{(m)}(1) \cr
&= (1-\la) \sum_{i\ge 1} {(\nabla g^{(m)})_i(1) \, \la^i \over 1-\la^i}. \cr
}$$
After differentiating $g_i(t)$ with respect to $t$ ($m$ times), then evaluating the result at $t=1$, 
and finally differencing with respect to $i$, the result is a bivariate polynomial 
$\calP(i, u)$ (with coefficients that are rational functions of $\la$) in the variables $i$ and $u = \la^i$.
Dividing this polynomial by $1-u = 1-\la^i$, we obtain $\calP(i,u) = \calQ(i,u)(1-u) + \calR(i)$, with quotient $\calQ(i,u)$ and remainder $\calR(i)$.
We then have
$$\Ex[N(N-1)\cdots(N-m+1)] = (1-\la)\sum_{i\ge 1} \calQ(i,\la^i) \, \la^i
+ (1-\la) \sum_{i\ge 1} {\calR(i) \, \la^i \over 1-\la^i}. \eqno(3.2)$$
The first sum in (3.2) can be expressed as a linear combination (with coefficients that are rational functions of $\la$) of sums of the form
$$S_{l,k}(\la) = \sum_{i\ge 1} i^l \, \la^{ki}.$$
These sums are themselves rational functions of $\la$:
$$S_{l,k}(\la) = {A_l(\la^k) \over (1-\la^k)^{l+1}},$$
where $A_l(x) = \sum_{0\le k\le l} a(l,k)\,x^k$ is the $l$-th Eulerian polynomial and the
$a(l,k)$ are the Eulerian numbers, with generating function $\sum_{l,k\ge 0} a(l,k) \, z^k \, y^l / l!
= z(1-z) / (e^{y(1-z)} - z)$
(see Comtet [C, p.~245]).
The second sum in (3.2) can be expressed as a linear combination (again with coefficients that are rational functions of $\la$) of the sums $T_l(\la)$ given by (1.8), with asymptotic formulas given by (1.9) and (1.10), and with exact formulas given by (1.11) and (1.12).

We are now ready to obtain the moments of $W$.
Differentiating the identity (2.3) $m$ times, we obtain
$$M^{(m)}_W(s) = F_m\(M^{(1)}_X(s), \ldots, M^{(m)}_X(s); 
g^{(1)}(M_X(s)), \ldots, g^{(m)}(M_X(s))\),$$
where $F_m(x_1, \ldots, x_m; y_1,  \ldots, y_m)$ is the polynomial
$$F_m(x_1, \ldots, x_m; y_1,  \ldots, y_m) = m!
\sum_{1\le l\le m} x_l \sum_{e_1, e_2, \dots, e_m} 
\prod_{1\le k\le m} \left({y_k\over k!}\right)^{e_k},$$
and the inner sum is over all $e_1, e_2, \dots, e_m$ such that $e_1 + e_2 + \cdots + e_m = l$
and $e_1 + 2\,e_2 + \cdots + m\,e_m = m$
(see for example Comtet [C, p.~137]).
Evaluating at $s=0$ and using $M_X(0)=1$ and $M^{(l)}_X(0) = 1/(1+\la)^l$ for $l\ge 1$, we have
$$\eqalign{
M^{(m)}_W(0) 
&= F_m\left(g^{(1)}(1), \ldots, g^{(m)}(1); {1\over 1+\la},\ldots, {1\over(1+\la)^m}\right) \cr
&= {1\over (1+\la)^m} \, \sum_{1\le l\le m} g^{(l)}(1) \sum_{e_1, e_2, \dots, e_m} 
{m\choose e_1, e_2, \dots, e_m} \cr
&= {1\over (1+\la)^m} \, \sum_{1\le l\le m} g^{(l)}(1) \, \Stwo{m}{l}, \cr
}$$
where the $\Stwo{m}{l}$
are the Stirling numbers of the second kind, with the generating function
$\sum_{m\ge l\ge 0} \Stwo{m}{l} {y^l \, z^m / l!} = e^{y(e^z-1)}$
(see for example Comtet [C, pp.~206--207]). 
Thus we have
$$\eqalign{
\Ex[W^m]
&= M_W^{(m)}(0) \cr
&= {1\over (1+\la)^m} \, \sum_{1\le l\le m} \Stwo{m}{l} \, g^{(l)}(1), \cr
}$$
where the $g^{(l)}(1)$ are the factorial moments of $N$ determined in the preceding paragraph.

The asymptotic formulas given above for the moments of $W$ can be extended to complete asymptotic expansions, with error terms of the form $O\((1-\la)^R\)$ for any $R$.
Any rational function of $\la$ has a Laurent series around $\la=1$, which will serve an an asymptotic expansion as $\la\to 1$ as well.
Thus the only remaining problem is to find asymptotic expansions for the sums
$T_l(\la)$.
These expansions have been given by Eschenfeldt, Gross and Pippenger [E].
We have
$$T_0(\la)\sim {1\over h}\log {1\over 1-\la} + {\ga\over h} 
+ \sum_{r\ge 0} {(-1)^r \, B_{r+1} \(B_{r+1} - (-1)^{r+1}\) \, h^r \over (r+1)\,(r+1)!},$$
where $\ga = 0.5772\ldots$ is Euler's constant, $B_r$ is the $r$-th Bernoulli number, defined by $t/(e^t - 1) = \sum_{r\ge 0} B_r \, t^r/r!$ (see for example Roman [R3, p.~94], and $h = -\log\la$ has the expansion
$$\eqalign{
h
&= -\log\(1 - (1-\la)\) \cr
&= \sum_{r\ge 1} {(1-\la)^r \over r}, \cr
}$$
so that its reciprocal has the expansion
$$\eqalign{
{1\over h}
&=  {1\over 1-\la}\, \sum_{r\ge 0} {(-1)^r \, C_r \, (1-\la)^r \over r!}, \cr
}$$
where $C_r$ is the $r$-th Bernoulli number of the second kind, defined by
$t/\log(1+t) = \sum_{r\ge 0} C_r \, t^r/r!$ (see for example Roman [R3, p.~116]).
(These numbers are also called the Cauchy numbers of the first kind, and are given by
$C_r = \int_0^1 x(x-1)\cdots(x-r+1)\,dx$; see for example Comtet [C, pp.~293--294].)
For $\l\ge 1$, we have
$$T_l(\la)\sim \sum_{r\ge 0} {(-1)^{r+l-1} \, B_r \, B_{r+l} \, h^{r-l} \over r! \, (r+l)}.$$
We note that, if $l$ is odd, then this expansion has only finitely many terms (because $B_r=0$ for odd $r\ge 3$).
\sk

\heading{5. Acknowledgment}

The research reported here was supported
by Grant CCF  0917026 from the National Science Foundation.
\vfill\eject

\heading{6. References}

\refbook C; L. Comtet;
Advanced Combinatorics:
The Art of Finite and Infinite Expansion;
D.~Reidel Publishing Co., Dortrecht, 1974. 

\item{[E]} 
P. Eschenfeldt, B. Gross and N. Pippenger,
``Stochastic Service Systems, Random Interval Graphs and Search Algorithms'',
{\tt arXiv:1107.4113v2 [math.PR]}.

\refbook G; G. Gasper and M. Rahman;
Basic Hypergeometric Series;
Cambridge University Press, Cambridge, 1990.

\refbook H; G. H. Hardy and E. M. Wright;
Introduction to the Theory of Numbers (5th edition);
Clarendon Press, Oxford, 1979.

\ref K; J. F. C. Kingman;
``The Effect of Queue Discipline on Waiting Time Variance'';
Math.\ Proc.\ Cambridge Phil.\ Soc.; 58:1 (1962) 163--164.

\ref L; J. D. C. Little;
``A Proof for the Queuing Formula: $L = \la W$'';
Oper.\ Res.; 9 (1961) 383--387.

\ref R1; J. Riordan;
``Delay Curves for Calls Served at Random'';
Bell System Technical Journal; 32 (1953) 100-119.

\refbook R2; J. Riordan;
Stochastic Service Systems;
John Wiley and Sons, New York, 1962.

\refbook R2; S. Roman;
The Umbral Calculus;
Academic Press, New York, 1984.

\ref T; D. G. Tambouratzis;
``On a Property of the Variance of the Waiting Time of a Queue'';
J. Appl. Prob.; 5:3 (1968) 702--703.

\ref V1; \'{E}. Vaulot;
``D\'{e}lais d'attente des appels t\'{e}l\'{e}phoniques trait\'{e}s au hasard'';
Comptes Rendus Acad.\ Sci.\ Paris; 222 (1946) 268--269.

\ref V2; \'{E}. Vaulot;
``D\'{e}lais d'attente des appels t\'{e}l\'{e}phoniques dans l'ordre inverse de leurs arriv\'{e}e'';
Comptes Rendus Acad.\ Sci.\ Paris; 238 (1954) 1188--1189.

\ref W1; B. Weesakul;
``The Random Walk between a Reflecting and an Absorbing Barrier'';
Ann.\ Math.\ Statist.; 32:3 (1961) 765--769.

\refbook W2; E. T. Whittaker and G. N. Watson;
A Course of Modern Analysis (4th edition);
Cambridge University Press, London, 1927.

\bye